\documentclass{amsart}
\usepackage{amsfonts}
\usepackage{amsmath,amssymb}
\usepackage{amsthm}
\usepackage{amscd}
\usepackage{graphics}
\usepackage{graphicx}

\usepackage[pagewise]{lineno}
\theoremstyle{definition}{
\newtheorem{Def}{{\rm Definition}}
\newtheorem{Ex}{{\rm Example}}
\newtheorem{Rem}{{\rm Remark}}

\newtheorem*{MainProb}{Main Problem}
}
\theoremstyle{plain}
{

\newtheorem{Thm}{Theorem}
\newtheorem{MainThm}{Main Theorem}

}

\begin{document}
\title[Reeb graphs of Morse functions with prescribed preimages]{Realization problems of graphs as Reeb graphs of Morse functions with prescribed preimages}
\author{Naoki Kitazawa}
\address{Institute of Mathematics for Industry, Kyushu University, 744 Motooka, Nishi-ku Fukuoka 819-0395, Japan\\
 TEL (Office): +81-92-802-4402 \\
 FAX (Office): +81-92-802-4405 \\
}
\email{n-kitazawa@imi.kyushu-u.ac.jp}
\urladdr{https://naokikitazawa.github.io/NaokiKitazawa.html}
\maketitle
\begin{abstract}
The present paper is on a new result on so-called realization problems of graphs as the {\it Reeb graphs} of Morse functions with prescribed preimages. The {\it Reeb graph} of a smooth function is
a graph obtained by identifying two points in the manifold of the domain if and only if they are in a same connected component of some preimage. In cases such as cases where functions are {\it Morse-Bott} functions on closed manifolds we obtain such graphs. They represent the manifolds of the domains compactly. They are not only fundamental tools in geometry using differentiable maps as fundamental tools, but also in applications of mathematics such as visualizations. Realization problems are also important. Constructing explicit functions with prescribed preimages was considered first by the author essentially and we generalize some of our related previous results for cases where manifolds of the domains are $3$-dimensional.
  
\end{abstract}

\section{Introduction}
\label{sec:1}
We introduce fundamental notation and terminologies on (differentiable) manifolds and (differentiable) maps between them.
For a differentiable manifold $X$, $T_xX$ denotes the tangent space at $x$. 
A point $x \in X$ in the manifold of the domain of a differentiable map $c:X \rightarrow Y$ is said to be {\it singular} if the rank of the differential $dc_x:T_xX \rightarrow T_{c(x)}Y$ is smaller than $\min \{\dim X,\dim Y\}$. 
We define the {\it singular set} of $c$ as the set of all singular points of $c$, denoted by $S(c)$.
Points in the image $c(S(c))$ are called {\it singular values} of $c$.
{\it Regular values} of the map are points in the manifold of the target which are not singular values of $c$.

${\mathbb{R}}^k$ denotes the $k$-dimensional Euclidean space for $k \geq 1$, endowed with the natural differentiable structure and the standard metric. 
${\mathbb{R}}^1$ is also denoted by $\mathbb{R}$ as usual. For $x \in {\mathbb{R}}^k$, $||x|| \geq 0$ denotes the distance between the origin $0 \in {\mathbb{R}}^k$ and $x$.
$S^k:=\{x \mid x \in {\mathbb{R}}^{k+1}, ||x||=1.\}$ denotes the $k$-dimensional unit sphere for $k \geq 0$. This is a 
compact smooth submanifold with no boundary (in ${\mathbb{R}}^{k+1}$) and connected for $k>0$. $D^k:=\{x \mid x \in {\mathbb{R}}^{k}, ||x|| \leq 1.\}$ denotes the $k$-dimensional unit disk for $k \geq 1$. This is a connected smooth submanifold in ${\mathbb{R}}^k$. 
\begin{Def}
\label{def:1}
A {\it Morse} function $c:X \rightarrow \mathbb{R}$ is a smooth function on a smooth manifold $X$ satisfying the following two  conditions.
\begin{enumerate}
\item Singular points are always in the interior of the manifold. 
\item At each singular point $x$ of which it has the form $$(x_1,\cdots,x_m) \mapsto {\Sigma}_{j=1}^{m-i(p)} {x_j}^2-{\Sigma}_{j=1}^{i(p)} {x_{m-i(p)+j}}^2+c(x)$$
for an integer $0 \leq i(p) \leq m$, which is shown to be uniquely chosen, and for suitable coordinates.
\end{enumerate}
A {\it Morse-Bott} function $c:X \rightarrow \mathbb{R}$ is a smooth function represented as the composition of a smooth submersion with a Morse function at each singular point.
\end{Def}
\subsection{Reeb spaces and Reeb graphs}
From a differentiable map $c:X \rightarrow Y$, we can define an equivalence relation ${\sim}_c$ on $X$ as follows: $x_1 {\sim}_{c} x_2$ holds if and only if they are in a same connected component of $c^{-1}(y)$ for some point $y$.

\begin{Def}
The quotient space $W_c:=X/{\sim}_c$ is defined as the {\it Reeb space} of $c$.
\end{Def}
For the Reeb space of $c$, $q_c:X \rightarrow W_c$ denotes the quotient space. We can define a map which is denoted by $\bar{c}$ uniquely so that we have the relation $c=\bar{c} \circ q_c$. 

Reeb spaces are in considerable cases graphs. For such a function $c$, the vertex set of $W_c$ is the set of all points $p \in W_c$ giving the preimages ${q_c}^{-1}(p)$ containing at least one singular point of the map $c$. $W_c$ is called the {\it Reeb graph} of $c$.
For example, for smooth functions on compact manifolds with finitely many singular values, Reeb spaces are homeomorphic to graphs and if they are on closed manifolds, such graphs are obtained (\cite{saeki}). Note that a graph is naturally topologized and (PL) homeomorphic to a $1$-dimensional polyhedron.
The present paper concentrates on such smooth functions on closed manifolds.

Reeb spaces represent the manifolds of the domains compactly. In specific cases, they inherit topological information such as homology groups and cohomology rings much, shown explicitly in \cite{saekisuzuoka} and papers and preprints \cite{kitazawa1,kitazawa2,kitazawa3,kitazawa6,kitazawa7,kitazawa8} of the author, for example. So-called {\it fold} maps such that preimages of regular values are disjoint unions of spheres are studied for example. In our present paper, we do not need to know these studies and results. {\it Fold} maps are, in short, smooth maps from smooth manifolds with no boundaries into manifolds with no boundaries represented as the product maps of Morse functions and identity maps on the interiors of copies of a unit disk locally (for suitable coordinates). The Reeb spaces of fold maps are shown to be polyhedra whose dimensions are same as those of the manifolds of the targets: see \cite{shiota} for general related theory and see also \cite{kobayashisaeki} for cases where the manifolds of the targets are surfaces, for example. We do not concentrate on {\it fold} maps in the present study.

For Reeb graphs and Reeb spaces, there exist various papers. \cite{reeb} seems to be one of pioneering papers. Other papers will be presented in the present paper.
\subsection{Some notions on graphs which are topologized canonically.}

\begin{Def}
An {\it isomorphism} between two graphs $K_1$ and $K_2$ is a (PL or piecewise smooth) homeomorphism from $K_1$ to $K_2$ mapping the vertex set of $K_1$ onto the vertex set of $K_2$ where they are topologized naturally.
\end{Def}
\begin{Def}
A continuous real-valued function $g$ on a graph $K$ is said to be a {\it good} function if it is injective on each edge.
\end{Def}

Hereafter, we only consider finite (and connected) graphs. A finite graph has such a function if and only if it has no loops as edges.

\subsection{Our main theorem.}
The present paper studies another important problem as follows.
\begin{MainProb}
For a finite and connected graph with at least one edge which may be a (so-called) multigraph, can we construct a smooth function on a (compact and connected) manifold (satisfying some good conditions) which induces the Reeb graph isomorphic to the graph and preimages of regular values for which are as prescribed? More explicitly, can we construct such functions as Morse functions?
\end{MainProb}
This is so-called realization problems of graphs as the Reeb graphs of smooth functions of suitable classes. \cite{sharko} is a pioneering paper on this. This concerns smooth functions on closed surfaces whose Reeb graphs are isomorphic to given graphs. Other related papers are presented in Remark \ref{rem:2}, in the last, and in other scenes in our paper, for example.

For the family of (all) smooth manifolds, we can induce the following equivalence relation: two smooth manifolds are equivalent if and only if they are diffeomorphic each other. 
The equivalence class which a smooth manifold belongs to is said to be the {\it diffeomorphism type} for the manifold.
A ({\it most fundamental}) {\it handlebody} is a smooth manifold of a suitable class of compact, connected and smooth manifolds whose dimensions are greater than $1$. The notion of a ({\it most fundamental}) {\it handlebody} will be presented in the next section. In short, this is a compact manifold obtained by attaching finitely many manifolds which are also so-called {\it handles} along the boundary of a copy of a unit disk disjointly.
In the present paper, we show the following result as Main Theorem \ref{mthm:1}.
\begin{MainThm}
\label{mthm:1}
Let $K$ be a finite and connected graph which has at least one edge and no loops. Let there exist a good function $g$ on $K$. Suppose that a diffeomorphism type for some manifold is assigned to each edge by a map $r_K$ on the edge set $E$. Let $m>1$ be an integer. Suppose the following conditions.
\begin{itemize}
\item At each edge of $K$ containing a vertex where $g$ has a local extremum, the value of $r_K$ is the diffeomorphism type for $S^{m-1}$.
\item Each vertex where $g$ has a local extremum is of degree $1$.
\end{itemize}
Assume also that the values of $r_K$ are always diffeomorphism types for closed and connected manifolds diffeomorphic to the boundaries of some $m$-dimensional most fundamental handlebodies. Then there exists an $m$-dimensional closed and connected manifold $M$ and a Morse function $f:M \rightarrow \mathbb{R}$ enjoying the following three properties.
\begin{enumerate}
\item \label{mthm:1.1} The Reeb graph $W_f$ of $f$ is isomorphic to $K${\rm :} we can take a suitable isomorphism $\phi:W_f \rightarrow K$ compatible with the remaining properties in these three properties.
\item \label{mthm:1.2} If we consider the natural quotient map $q_f:M \rightarrow W_f$ and for each point $\phi(p) \in K$ {\rm (}$p \in W_f${\rm )} that is not a vertex and that is in an edge $e$, then the diffeomorphism type for
 the preimage ${q_f}^{-1}(p)$ is $r_K(e)$.
\item \label{mthm:1.3} For each point $p \in M$ mapped by $q_f$ to a vertex $v_p:=q_f(p) \in W_f$, $f(p)=g \circ \phi(v_p)$.

\end{enumerate}
\end{MainThm}
The next section explains about handlebodies. The third section presents one of our motivating results on smooth functions on $3$-dimensional closed manifolds.
We prove Main Theorems in the fourth section. 
Roughly speaking, methods are similar to the methods in the proof of Theorem \ref{thm:1} or the main theorem of \cite{kitazawa4} (and methods in the proofs of main theorems of  \cite{kitazawa5, kitazawa9}). However, related differential topological arguments on manifolds belonging to certain classes we need are mutually different.
Main Theorem \ref{mthm:2} applies Main Theorem \ref{mthm:1} where the dimensions are $m=4$. \\
\ \\
{\bf Acknowledgement, grants and data.} \\
\thanks{The author was a member of the projects JSPS KAKENHI Grant Number JP17H06128 and JP22K18267. Principal Investigators are all Osamu Saeki. This work was also supported by the project. He is not a member of JSPS KAKENHI Grant Number JP23H05437 formally. However members related to this project support our studies. Principal investigator is also Osamu Saeki. The speaker is a researcher at Osaka Central
		Advanced Mathematical Institute (OCAMI researcher). He is not employed there. This is for our studies and helps our studies. 

	The author would like to thank Osamu Saeki again for private discussions on \cite{saeki} with \cite{kitazawa4}. These discussions continue to motivate the author to study the present study and related studies further. In addition,
the author would like to thank him for explicit positive comments on Main Theorem \ref{mthm:2}. He has said that the result must be important in (differential) topological theory of $3$-dimensional and $4$-dimensional manifolds for example. 
The author would like to thank Takahiro Yamamoto for asking several questions on the talk of the author related to the present study in the conference "Singularity theory of smooth maps and its applications" (http://www.math.kobe-u.ac.jp/HOME/saji/math/conf2021/indexe.html). This conference is supported by the Research Institute for Mathematical Sciences, an International Joint Usage/Research Center located in Kyoto University. He asked the author about the orientability of the manifolds $M$ in Theorem \ref{thm:1} and Main Theorem \ref{mthm:2} for example. The author would like to thank him again for additional informal discussions on the present paper with \cite{kitazawa4, kitazawa5} and the talk and his constant encouragement. He kindly presented some new notions on Reeb graphs such as the {\it decorated Reeb graphs} of smooth functions, respecting smooth functions on manifolds with non-empty boundaries, as a related and advanced topic. He also told the author that this is closely related to \cite{saekiyamamoto} for example.
 
All data supporting the present study are in the present paper.}
\section{Handlebodies.}
We review several notions on {\it handlebodies}.

A {\it k-handle} is a smooth manifold diffeomorphic to a smooth manifold of the form $D^k \times D^{m-k}$ ($0 \leq k \leq m$) as a cornered manifold where $D^0$ denotes a one-point set. We only consider the case $1 \leq k \leq m-1$ essentially.

A {\it linear} bundle is a bundle whose fiber is diffeomorphic to a Euclidean space, a unit sphere or a unit disk and whose structure group is linear and acts in a canonical and linear way on the fiber.

A {\it normal bundle} is a linear bundle over a smoothly embedded closed manifold in a smooth closed manifold and defined in a natural way based on elementary arguments on differential topology. The fiber is the Euclidean space whose dimension is same as the codimension of the submanifold. The total space is regarded as a smooth submanifold whose codimension is $0$. By restrcitng the fiber to the unit disk of the same dimension, we have a smooth compact submanifold whose codimension is $0$ naturally.  We can define such notions in suitably generalized scenes. We do not explain about the notions precisely and we also expect that we have some related knowledge.

As a fundamental fact on differential topology of manifolds, it is well-known that we can always eliminate corners for smooth manifolds in certain canonical ways and obtain smooth manifolds with no corners. Furthermore, the resulting diffeomorphism types are unique.

We introduce an operation of constructing compact smooth manifolds of good classes.
Take an integer $m>1$.
\begin{enumerate}
\item Set ${H_0}^m$ as a non-empty space which is also an $m$-dimensional smooth compact manifold.
\item We do either of the following two.
\begin{enumerate}
\item Set $k=0$ and take two distinct points in $\partial {H_0}^m$. We choose a small closed tubular neighborhood $N_0$ for the discrete set consisting of the two points. This is always seen as a normal bundle which is trivial and diffeomorphic to the disjoint union $D^{m-1} \sqcup D^{m-1}$.
\item Set a positive integer $k<m-1$ and take a smooth submanifold diffeomorphic to $S^k$ whose suitable small closed tubular neighborhood $N_0$ is regarded as a trivial linear bundle over this copy of $S^k$ seen as a normal bundle.
\end{enumerate}
\item Take a ($k+1$)-handle, identified with $D^{k+1} \times D^{m-k-1}$ via a diffeomorphism between the cornered manifolds.  
\item Attach the ($k+1$)-handle $D^{k+1} \times D^{m-k-1}$ just before to $\partial {H_0}^m$ via a bundle isomorphism between the trivial linear bundle $\partial D^{k+1} \times D^{m-k-1} \subset D^{k+1} \times D^{m-k-1}$ whose fiber is diffeomorphic to $D^{m-k-1}$ and the small closed tubular neighborhood $N_0$.
\item We eliminate the corner of the resulting manifold in the canonical way and set the resulting manifold as ${H_1}^m$.
\item We do a finite iteration of this procedure and obtain $m$-dimensional smooth manifolds with no corners ${H_{l_0}}^m$ for $l_0 \geq 1$ where $l_0$ is an integer. At each step we first choose an $m$-dimensional compact and smooth manifold ${H_{l_0-1}}^m$ instead of "${H_0}^m$" and the closed tubular neighborhood $N_{l_0-1}$ instead of "$N_0$".
\end{enumerate}
\begin{Def}
In this operation, let ${H_0}^m$ be diffeomorphic to the unit disk. We call the manifold ${H_{l_0}}^m$ a {\it handlebody} for $l_0 \geq 1$. In addition, suppose the following two.
\begin{itemize}
\item ${H_0}^m$ is diffeomorphic to the unit disk.
\item The small closed tubular neighborhoods denoted by the notation $N_{l-1}$ are always in the original boundary $\partial {H_0}^m \subset {H_0}^m$ and mutually disjoint for $1 \leq l \leq l_0$.
\end{itemize}
In this case, we call the $m$-dimensional manifold ${H_{l_0}}^m$ a {\it most fundamental} handlebody.
\end{Def}

\begin{Ex}
The boundary of a copy of the $m$-dimensional unit disk $D^m$ and equivalently, a copy of the unit sphere $S^{m-1}$, are also the boundaries of some most fundamental handlebodies for $m \geq 2$. We attach two handles to the boundary of a copy ${H_0}^m$ of the $m$-dimensional unit disk simultaneously to obtain this. More precisely, we attach a $1$-handle and an ($m-1$)-handle to obtain the desired handlebody by defining the closed tubular neighborhoods $N_0$ and $N_1$ before satisfying the following two conditions.
\begin{enumerate}
\item $N_0$ is the union of a suitable two copies of $D^{m-1}$ smoothly and disjointly embedded in $\partial {H_0}^m$. 
\item $N_1$ is the closed tubular neighborhood of a copy of $S^{m-2}$ smoothly embedded in $\partial {H_0}^m$ and two connected components of $N_0$ are in distinct connected components of $ \partial {H_0}^m-N_1$. Here we consider each $H_j$ as a smooth compact submanifold $H_j \subset H_{j+1}$ in a natural way.
\end{enumerate}
We have ${H_2}^m$ and $\partial {H_2}^m$ is diffeomorphic to $\partial D^m=S^{m-1}$.
\end{Ex}

\section{A previous motivating result.}
The following theorem is previously obtained and motivates us to study the present problem. Main Theorem \ref{mthm:1} can be regarded as an extension of Theorem \ref{thm:1} except the last part with the orientability of $M$.
\begin{Thm}[\cite{kitazawa4}]
	\label{thm:1}
	Let $K$ be a finite and connected graph which has at least one edge and no loops. Let there exist a good function $g$ on $K$. Suppose that a non-negative integer is assigned to each edge by a map $r_K$ on the edge set $E$. These two functions satisfy the following conditions.
	\begin{itemize}
		\item At each edge of $K$ containing a vertex where $g$ has a local extremum, the value of $r_K$ is $0$.
		\item Each vertex where $g$ has a local extremum is of degree $1$.
	\end{itemize}
	Then there exist a $3$-dimensional closed, connected and orientable manifold $M$ and a Morse function $f:M \rightarrow \mathbb{R}$ enjoying the following three properties.
	\begin{enumerate}
		\item The Reeb graph $W_f$ of $f$ is isomorphic to $K$ and we can take a suitable isomorphism $\phi:W_f \rightarrow K$ compatible with the remaining properties of these three properties.
		\item For each point $\phi(p) \in K$ {\rm (}$p \in W_f${\rm )} in the interior of an edge $e$, the preimage ${q_f}^{-1}(p)$ is a closed, connected, and orientable surface of genus $r_K(e)$.
		\item For each point $p \in M$ mapped by $q_f$ to a vertex $v_p:=q_f(p) \in W_f$, $f(p)=g \circ \phi(v_p)$.
	\end{enumerate}
	Furthermore, if we drop the two conditions on $g$, $r_K$ and the graph $K$, then there exist a $3$-dimensional closed, connected and orientable manifold $M$ and a smooth function $f$ on $M$ enjoying the three properties just before and the following two properties are also enjoyed.	\begin{enumerate}
		\item At each singular point where $f$ does not have a local extremum, it is locally a Morse function.
		\item At each singular point where $f$ has a local extremum, except finitely many such singular points, it is locally a Morse-Bott function. 
	\end{enumerate}
\end{Thm}

Our preprint \cite{kitazawa9} is an extension of this to a non-orientable case.
\section{On Main Theorems.}
We present some fundamental arguments on handlebodies, Morse functions and differential topology of manifolds. Any smooth, compact and connected manifold always has a Morse function such that at distinct singular points, the values are always distinct.

Choose an arbitrary disjoint ordered pair of non-empty subsets represented as disjoint unions of connected components of the boundary such that their union is the original boundary. We always have such a Morse function such that the preimage of the minimum and that of the maximum coincide with the former submanifold and the latter submanifold, respectively.

If either of the subsets is empty instead here, then we have a Morse function with similar properties such that either the preimage of the minimum or the maximum coincides with the boundary.


It is also well-known that Morse functions enjoying the presented properties exist densely in the space of all smooth functions topologized suitably. Related to this, for fundamental theory on singularity theory and differential topological properties of Morse functions (and fold maps), see \cite{golubitskyguillemin} for example.

Let us review a correspondence between singular points of Morse functions and handles. This plays essential roles in our construction of functions including our proof of Main Theorem \ref{mthm:1}.

From a Morse function as presented here on a compact manifold $E$ whose boundary $\partial E$ is non-empty and represented as the disjoint union $F_1 \sqcup F_2$ where $F_j$ is not empty and the disjoint union of finitely many connected components of $\partial E$. We have an explicit decomposition of $E$ into ($k$-)handles in the construction of the $m$-dimensional manifold ${H_{l_0}}^m$ before. We can set $F_1$, $F_2$ and ${H_0}^m$ so that the following conditions hold.
\begin{itemize}
\item ${H_0}^m$ is a (small) collar neighborhood $F_1 \times [-1,0]$ of $F_1$, identified with $F_1 \times \{-1\} \subset F_1 \times [-1,0]$ (in a canonical way).
\item $F_1 \times \{-1\}$ is the preimage of the minimum and that $F_2$ is the preimage of the maximum.
\end{itemize}
In this situation, each $k$-handle naturally corresponds to a singular point $p$ such that $i(p)=k$ in Definition \ref{def:1} and for the $j_1$-th handle and the $j_2$-th handle with $j_1<j_2$, the corresponding (singular) value at the singular point corresponding to the $j_2$-th singular point is greater than the corresponding value at the singular point corresponding to the $j_1$-th singular point.

Conversely, from a manifold ${H_0}^m$ diffeomorphic to $F_1 \times [-1,0]$ and such a system of handles for an $m$-dimensional compact and smooth manifold $E:={H_{l_0}}^m$ satisfying $\partial E=(F_1 \times \{-1\}) \sqcup F_2$, we have a Morse function enjoying the properties where $F_1$ and $F_2$ are non-empty and closed manifolds which may not be connected.

As a specific case, we consider a compact manifold $E$ such that the following three hold where we abuse notation in the construction of handlebodies and more general compact manifolds.
\begin{itemize}
\item The boundary $\partial E$ is represented as the disjoint union $F_1 \sqcup F_2$ where $F_j$ is not empty and the disjoint union of finitely many connected components of $\partial E$. 
\item Small closed tubular neighborhoods denoted by $N_j$ ($0 \leq j \leq l_0-1$) are always in the boundary $\partial {H_0}^m$ of ${H_0}^m \subset E$ and mutually disjoint as in the definition of a most fundamental handlebody in Definition \ref{def:1}.
\item ${H_0}^m$ is a (small) collar neighborhood $F_1 \times [-1,0]$ of $F_1$, identified with $F_1 \times \{-1\}$ in a canonical way.
\end{itemize}
In this case we have a Morse function with exactly one singular value instead where there may exist singular points at least $2$. The singular set is a discrete subset of course and a finite subset. Conversely, from such a Morse function, we have a system of handles as presented here.

Related theory is explained in \cite{milnor} for example. This is also a key ingredient in \cite{michalak}. 
In the paper Michalak has constructed Morse functions such that connected components of preimages containing no singular points (of the functions) are diffeomorphic to unit spheres. Earlier versions of \cite{kitazawa4} adopt expositions via this theory. 
\cite{kitazawa9} is, as presented shortly in the third section, an extension of \cite{kitazawa4}. This finds certain sufficient conditions to construct smooth functions which are locally Morse around singular points and whose preimages having no singular points (of the functions) may be non-orientable closed surfaces. 
See also Remark \ref{rem:2}, presented in the end of our present paper. We do not need to understand related arguments precisely. We discuss important arguments in our present paper in a self-contained way. 

We go to our proof of Main Theorem \ref{mthm:1}. We also abuse our notation in the construction of handlebodies and more general compact manifolds, presented in the second section and this section. 

\begin{proof}[A proof of Main Theorem \ref{mthm:1}.]
As the original proof of Theorem \ref{thm:1}, we construct a local Morse function around each vertex, respecting the three properties (\ref{mthm:1.1})--(\ref{mthm:1.3}) in Main Theorem \ref{mthm:1} and glue them together.
We construct the local functions. After that, we also glue them together to complete the proof.

Around each vertex $v$ where $g$ has a local extremum, we construct a Morse function on a copy of the unit disk obtained by considering the natural height. This is a so-called {\it height} function. This function $f_{v,m}:D^{m} \rightarrow \mathbb{R}$ is defined by $f_{v,m}(x):=\pm||x||^2+g(v)$ and $0$ is mapped to $v$ by the composition of $q_{f_{v,m}}$ with a suitable PL homeomorphism ${\phi}_v$ from the Reeb space onto a small regular neighborhood of $v$ in $K$.

Around each vertex $v$ where $g$ does not have a local extremum, we construct a desired local Morse function. This presents a new ingredient.

Set ${H_0}^m$ as a copy of the $m$-dimensional unit disk $D^m$ in the construction of a handlebody. We can attach $l \geq 0$ ($m-1$)-handles and additional suitable handles simultaneously as explained in the definition of a most fundamental handlebody to obtain a most fundamental handlebody whose boundary consists of exactly $l+1$ connected components so that the following conditions hold.
\begin{itemize}
\item We attach at least one handle here.
\item The boundary of the union of ${H_0}^m$ and the first $l$ ($m-1$)-handles is the disjoint union of $l+1$ copies of $S^{m-1}$ (after the corner is eliminated in the canonical way).
\item Each of the remaining handles is attached to (the interior of) the intersection of the original boundary $\partial {H_0}^m$ and one of the $l+1$ copies of $S^{m-1}$ before. Here we consider the natural identification.
\end{itemize}
The $l+1$ connected components of the resulting boundary are the boundaries of some most fundamental handlebodies by the definition. 

For any ($m-1$)-dimensional manifold $X^{m-1}$ diffeomorphic to the disjoint union of $l+1$ connected manifolds which are the boundaries of some most fundamental handlebodies, we can attach the handles to the boundary $\partial {H_0}^m \subset {H_0}^m$ of the copy of the unit disk $D^m$ so that the boundary of the resulting most fundamental handlebody is diffeomorophic to $X^{m-1}$ before. We can know this immediately from the construction and the definition of a most fundamental handlebody.

We can remove a smoothly embedded copy of the $m$-dimensional unit disk in the interior of the original copy ${H_0}^m$ of the unit disk $D^m$ (smoothly embedded) in a most fundamental handlebody obtained in such a way. In other words, we can change the manifold ${H_0}^m$ to a one diffeomorphic to $S^{m-1} \times [-1,0]$. 

We prepare two compact and connected manifolds obtained in such ways. Let $E_1$ and $E_2$ denote them. Let ${\partial}_0 {H_1} \subset \partial E_1$ and ${\partial}_0 {H_2} \subset \partial E_2$ denote the connected components of the boundaries originally copies of the unit disks in the interior of ${H_0}^m$ are attached to. In other words, $E_i$ is a manifold obtained by removing the interior of a copy of $D^{m}$ smoothly embedded in the interior of the submanifold ${H_0}^m$ in a most fundamental handlebody and ${\partial}_0 {H_i} \subset \partial E_i$ is diffeomorphic to $S^{m-1}$.

Let ${X_1}^{m-1} \subset \partial E_1$ and ${X_2}^{m-1} \subset \partial E_2$ denote the disjoint unions of the remaining connected components of the boundaries. By the exposition on the relationship between singular points of Morse functions and handles (with some additional arguments), we easily have a Morse function with exactly one singular value such that the preimage of the minimum is ${\partial}_0 H_1$, that the preimage of the maximum is ${X_{1}}^{m-1}$ and that the singular value is $g(v)$. By the same reason, we also have another Morse function with exactly one singular value such that the preimage of the maximum is ${\partial}_0 H_2$, that the preimage of the minimum is ${X_{2}}^{m-1}$ and that the singular value is $g(v)$. We can obtain these functions so that the images are $[m_v,M_v]$ and contain the unique singular values $g(v)$ in the interior $(m_v,M_v)$. We can choose a trivial linear bundle over $[m_v,M_v]$ whose fiber is diffeomorphic to $D^{m-1}$ apart from the singular set and remove the interior of the total space of this trivial bundle for each Morse function. We can do this thanks to the structures of the Morse functions. After that we can glue the resulting functions on the resulting manifolds preserving the values at each point in the manifolds of the domains. In this way, we can obtain a local Morse function $f_{v}$ enjoying the following properties.
\begin{itemize}
\item The manifold of the domain is compact and connected and its boundary is diffeomorphic to ${X_1}^{m-1} \sqcup {X_2}^{m-1}$.
\item The image is $[m_v,M_v]$.
\item The preimage of the maximum is diffeomorphic to ${X_1}^{m-1}$. The preimage of the minimun is diffeomorphic to ${X_2}^{m-2}$.
\item There exists exactly one singular value and the singular value is $g(v)$.
\item The Reeb space is homeomorphic to a graph. This is also regarded as a specific case of a main theorem of \cite{saeki}.
\item The preimage of the singular value $g(v)$ is connected ($q_{f_{v}}(p)$ in the next property is uniquely determined). 
\item There exists a suitable PL homeomorphism ${\phi}_v$ from the Reeb space onto a small regular neighborhood of $v$ in $K$ such that the composition of $q_{f_{v}}$ with ${\phi}_v$ maps $p$ to $v$ for any singular point $p$.
\end{itemize}.

Note that the pair $({X_1}^{m-1},{X_2}^{m-1})$ can be also chosen as an arbitrary pair of two non-empty sets represented as some disjoint unions of finitely many ($m-1$)-dimensional closed and connected manifolds diffeomorphic to the boundaries of some $m$-dimensional most fundamental handlebodies.

The conditions on $g$, $r_K$ and the graph $K$ imply that this completes the proof with the arguments in the original proof of Theorem \ref{thm:1}.
We explain about this more precisely. Around each connected component $R \subset K$ of the complementary set of the interior of the small regular neighborhood of the vertex set $V$ of $K$ (obtained as the disjoint union of the small regular neighborhoods of the vertices chosen before), we can construct a desired local function $f_R$ as a trivial smooth bundle over $R$: $R$ is in the interior of an edge $e_R \subset K$ and PL homeomorphic to a closed interval. We can easily have a PL or piecewise smooth homeomorphism ${\phi}_R$ from the Reeb space of the function to $R$.
 
We have a desired function $f$ and a desired PL or piecewise smooth homeomorphism $\phi$ by gluing the local smooth functions and the local PL or piecewise smooth homeomorphisms together.
\end{proof}

We consider the case $m=4$ in Main Theorem \ref{mthm:1}.

\begin{MainThm}
\label{mthm:2}

Let $K$ be a finite and connected graph which has at least one edge and no loops. Let there exist a good function $g$ on $K$. Suppose that a diffeomorphism type for some manifold is assigned to each edge by a map $r_K$ on the edge set $E$. Suppose the following conditions.
\begin{itemize}
\item At each edge of $K$ containing a vertex where $g$ has a local extremum, the value of $r_K$ is the diffeomorphism type for $S^{3}$.
\item Each vertex where $g$ has a local extremum is of degree $1$.
\item The values of $r_K$ are always diffeomorphism types for $3$-dimensional closed, connected and orientable manifolds.
\end{itemize}
Then there exists a $4$-dimensional closed, connected and orientable manifold $M$ and a Morse function $f:M \rightarrow \mathbb{R}$ enjoying the three properties {\rm (\ref{mthm:1.1})}--{\rm (\ref{mthm:1.3})} in Main Theorem \ref{mthm:1}.
\end{MainThm}
\begin{proof}
We need topological theory on $3$-dimensional manifolds such as so-called ({\it integral}) {\it Dehn surgeries}. For systematic explanations on the theory, see \cite{hempel} for example.

By the theory of (integral) Dehn-surgeries on $3$-dimensional manifolds, every $3$-dimensional closed, connected and orientable manifold can be regarded as the boundary of some $4$-dimensional most fundamental handlebody which is also orientable. (Our proof of) Main Theorem \ref{mthm:1} and an argument on the orientability of the manifold $M$ of the domain in the last of the original proof of Theorem \ref{thm:1} complete the proof. More precisely, we can glue the local functions one after another preserving the orientablility. 
\end{proof}

In \cite{kitazawa4, kitazawa9}, connected components of preimages containing no singular points of the constructed smooth functions are closed and connected surfaces. Note that it may not be represented as the boundary of any most fundamental handlebody in \cite{kitazawa9}. The real projective plane is of important examples. Closed and connected surfaces which are the boundaries of some most fundamental handlebodies are surfaces represented as connected sums of copies of $S^1 \times S^1$ or the total space of a linear bundle over $S^1$ whose fiber is $S^1$ and which is not trivial. $S^1 \times S^1$ is for the $2$-dimensional torus. The total space of a linear bundle over $S^1$ whose fiber is $S^1$ and which is not trivial is for the Klein Bottle.
In the $3$-dimensional case, $3$-dimensional closed and connected manifolds can not be represented as connected sums of copies of $S^1 \times S^2$ or the total spaces of linear bundles over $S^1$ or $S^2$ which are not trivial in general.

\begin{Rem}
	\label{rem:1}
For theory on diffeomorphism types for handlebodies and attachments of handles, see papers by Wall (\cite{wall1,wall2,wall3,wall4,wall5,wall6}) and see also \cite{ranicki} for example.
\end{Rem}
\begin{Rem}
	\label{rem:2}
We close the present paper by introducing related studies among the others with results of the author.

\cite{michalak} has given an answer to Main Problem by constructing Morse functions such that preimages of regular values (of the functions) are disjoint unions of copies of unit spheres. \cite{michalak2} studies deformations of Morse functions and their Reeb graphs. They motivated the author to give a new answer as \cite{kitazawa4} or Theorem \ref{thm:1} and \cite{kitazawa9}. Our preprint \cite{kitazawa9} is, as presented, an extension to similar smooth functions on (some) $3$-dimensional manifolds which may be non-orientable. In the preprint, preimages of the functions may be closed non-orientable surfaces. As presented, explicit sufficient conditions for the extension are found and discussed.

\cite{saeki} is regarded as a paper motivated by \cite{kitazawa4}. This also studies more general cases where diffeomorphism types for closed (compact), connected and smooth manifolds assigned to edges are general satisfying only conditions on so-called cobordism relations on closed (or compact) smooth manifolds. On the other hand, this does not study functions with very explicit singularities such as Morse(-Bott) functions.

\cite{gelbukh,gelbukh2} are related studies where Morse-Bott functions on closed surfaces and general closed manifolds are considered. These studies respecting preimages of regular values are new important problems and regarded as variants of Main Problem.

Last, \cite{kitazawa9} is also a study on Main Problem. This is on a case where preimages may not be compact. This is new in non-proper smooth functions are considered.
\end{Rem}


\end{document}